\documentclass[11pt,a4paper]{article}
\usepackage[USenglish]{babel}
\usepackage[T1]{fontenc}
\usepackage[utf8]{inputenc}
\usepackage{csquotes}
\usepackage{enumitem}
\usepackage{appendix}
\usepackage[numbers,sort]{natbib}
\usepackage{doi}
\usepackage{hyperref}

\hypersetup{colorlinks=true,linkcolor=blue,bookmarksnumbered=true,bookmarksdepth=2,pdftitle=Representation\ Formula\ for\ Viscosity\ Solution\ to\ a\ PDE\ Problem\ involving\ Pucci's\ Extremal\ Operator,pdfauthor=Marco\ Pozza}

\usepackage{sectsty}
\subsubsectionfont{\normalfont\itshape}

\usepackage{amssymb}
\usepackage{amsmath}
\usepackage{mathtools}
\usepackage{amsthm}
\usepackage{mathrsfs}
\usepackage{dsfont}

\usepackage[noabbrev]{cleveref}

\theoremstyle{plain}
\newtheorem{theo}{Theorem}[section]
\crefname{theo}{theorem}{theorems}
\crefname{dpp}{dynamic programming principle}{}
\newtheorem{prop}[theo]{Proposition}

\newtheorem{lem}[theo]{Lemma}
\crefname{lem}{lemma}{lemmas}
\theoremstyle{definition}
\newtheorem{defin}[theo]{Definition}

\newtheorem{prob}[theo]{Problem}

\crefname{assum}{assumptions}{assumptions}
\theoremstyle{remark}
\newtheorem{rem}[theo]{Remark}

\DeclareMathOperator{\tr}{tr}
\DeclareMathOperator*{\esss}{ess\,sup}

\DeclareMathOperator{\eig}{eig}
\DeclareMathOperator{\interior}{int}

\title{Representation Formula for Viscosity Solution to a PDE Problem involving Pucci's Extremal Operator}
\author{Marco Pozza\thanks{Dipartimento di Matematica ``G. Castelnuovo'', Sapienza Università di Roma. Piazzale Aldo Moro 5, 00185 Roma, Italy. {\tt pozza@mat.uniroma1.it}}}
\date{}

\begin{document}
\maketitle

\begin{abstract}
We provide a representation formula for viscosity solutions to an elliptic Dirichlet problem involving Pucci's extremal operators. This is done through a dynamic programming principle derived from \cite{art:pathpeng}. The formula can be seen as a nonlinear extension of the Feynman--Kac formula.
\end{abstract}


\medskip


\section{Introduction}

The purpose of this paper is to provide a representation formula for the viscosity solution to the problem
\begin{equation}\label{eq:pprob}
\left\{
\begin{aligned}
&\dfrac12\mathcal P^+_{\lambda,\Lambda}\left(D^2u\right)+f(x)=0,&\qquad&x\in D,\\
&u(x)=g(x),&&x\in\partial D,
\end{aligned}
\right.
\end{equation}
where $D$ is a bounded domain satisfying an \emph{exterior cone condition}, $f$ and $g$ are two continuous function in $\overline D$, and $\mathcal P^+_{\lambda,\Lambda}$ is Pucci's maximal operator.

Let us recall here the definition of Pucci's maximal operator: given two parameters $\Lambda\ge\lambda>0$, for any symmetric $N\times N$ matrix $S$,
\[
\mathcal P^+_{\lambda,\Lambda}(S):=\Lambda\sum_{\lambda_i\ge0}\lambda_i+\lambda\sum_{\lambda_i<0}\lambda_i,
\]
where $\{\lambda_i\}_{i=1}^N$ are the eigenvalues of $S$, or equivalently
\[
\mathcal P^+_{\lambda,\Lambda}(S):=\max_{A\in\mathcal M_{\lambda,\Lambda}^N}\tr(AS),
\]
where $\mathcal M_{\lambda,\Lambda}^N$ is the set containing the $N\times N$ symmetric matrix such that their eigenvalues are in $[\lambda,\Lambda]$. Pucci's extremal operators represent an important prototype of fully nonlinear operators and as such they are the main subject of many papers, as in e.g. \cite{cutrileoni00,quaassirakov06,felmerquaas05}. These operators also have a central role in the study of the regularity of viscosity solutions to fully nonlinear second order PDEs, as can be seen in the monograph of Cabré and Caffarelli \cite{CaCa95}. For a comprehensive treatment of the theory of viscosity solution we refer to \cite{userguide}.

We will prove that the viscosity solution to \eqref{eq:pprob} is given by the formula
\begin{equation}\label{eq:vissol}
u(x):=\sup_{\sigma\in\mathcal A}\mathds E\left(g\left(X^x_{\sigma,\tau_\sigma^x}\right)+\int_0^{\tau_\sigma^x}f\left(X_{\sigma,t}^x\right)dt\right),
\end{equation}
where $\mathcal A$ is a set that will be defined later, $\{X_\sigma\}_{\sigma\in\mathcal A}$ is a collection of stochastic processes and $\tau_\sigma$ is the \emph{exit time} of $X_\sigma$ from the domain $D$. The connection with the well known \emph{Feynman--Kac formula} is evident, indeed this is to all intents and purposes a nonlinear extension of the Feynman--Kac formula.

The main tool to prove \eqref{eq:vissol} is a \emph{dynamic programming principle}, which will be illustrated later and is an adaptation of the one presented by Denis, Hu and Peng in \cite{art:pathpeng}. It is obtained using stochastic control techniques, in this regard we cite, among the others, \cite{hjlions2,hjlions1,bensoussa82,krylovbook,flemingrishel75,nisio15}. We point out that in \cite{art:pathpeng} this principle is used to give a representation formula to a Cauchy type problem involving a sublinear operator which generalize Pucci's maximal operators and is related to nonlinear expectation, a topic thoroughly studied by Peng in this and in many other articles, see \cite{book:peng} for a complete overview of this subject. We recall that in \cite{pozza19}, we have used the same techniques to prove a representation formula for solutions to a more general class of nonlinear parabolic PDEs.

Another representation formula for solutions of the problem \eqref{eq:pprob} has been proposed by Blanc, Manfredi and Rossi in \cite{BlaJMaDRo19} under the additional assumption that $D$ satisfies a uniform exterior sphere condition. They approximate the viscosity solution to \eqref{eq:pprob} using a family of functions $\{u_\varepsilon\}$, where, for any $\varepsilon>0$, $u_\varepsilon$ is given as a supremum over all the possible strategies of the expected outcome of a game. Roughly speaking, we can say that $u_\varepsilon$ tends to the function represented by formula \eqref{eq:vissol} as $\varepsilon$ tends to 0.

We conclude pointing out that our method can be extended, going along the same lines of \cite{pozza19}, to obtain representation formulas for viscosity solutions to more general Dirichlet problems.

\subsubsection*{Notation}
\pdfbookmark[2]{Notation}{Notation}

Here we fix the notation and some conventions that we will use later. We also recall some basic results. Assume that we have a filtered probability space $(\Omega,\mathcal F,\{\mathcal F_t\}_{t\in[0,\infty)},\mathds P)$, then
\begin{itemize}
\item $\mathcal F$ is a complete $\sigma$--algebra on $\Omega$;
\item the stochastic process $\{W_t\}_{t\in[0,\infty)}$ will denote the $N$ dimensional Brownian motion under $\mathds P$;
\item $\{\mathcal F_t\}_{t\in[0,\infty)}$ is the filtration defined by $\{W_t\}_{t\in[0,\infty)}$ satisfying the usual condition of completeness and right continuity;
\item $\{W_s^t\}_{s\in[t,\infty)}:=\{W_s-W_t\}_{s\in[t,\infty)}$ is a Brownian motion independent from $\{W_s\}_{s\in[0,t]}$ by the \emph{strong Markov property};
\item $\{\mathcal F^t_s\}_{s\in[t,\infty)}$ is the filtration generated by $\{W_s^t\}_{s\in[t,\infty)}$ which we assume satisfy the usual condition and is independent from $\mathcal F_t$;
\item we will say that a stochastic process $\{H_t\}_{t\in[0,\infty)}$ is adapted if $H_t$ is $\mathcal F_t$--measurable for any $t\in[0,\infty)$;
\item we will say that a stochastic process $\{H_t\}_{t\in[0,\infty)}$ is progressively measurable, or simply progressive, if, for any $T\in[0,\infty)$, the application that to any $(t,\omega)\in[0,T]\times\Omega$ associate $H_t(\omega)$ is $\mathcal B([0,T])\times\mathcal F_T$--measurable;
\item a function on $\mathds R$ is called \emph{cadlag} if is right continuous and has left limit everywhere;
\item a cadlag (in time) process is progressive if and only if is adapted;
\item $\mathcal F_\infty:=\sigma(\mathcal F_t|t\in[0,\infty))$;
\item if $A\in\mathds R^{N\times M}$ then $A^\dag$ will denote its transpose and $\eig_A$ its spectrum;
\item (Frobenius product) if $A,B\in\mathds R^{N\times M}$ then
\[
\langle A,B\rangle:=\tr\left(AB^\dag\right)=\sum\limits_{i=1}^N\sum\limits_{j=1}^MA_{i,j}B_{i,j};
\]
\item if $A\in\mathds R^{N\times M}$ then $|A|$ will denote the norm $\sqrt{\langle A,A\rangle}=\sqrt{\sum\limits_{i=1}^N\sum\limits_{j=1}^MA_{i,j}^2}$;
\item $\mathds S^N$ is the set containing all the symmetric matrix of $\mathds R^{N\times N}$ and $\mathds S^N_+$ is its the subset containing all the positive definite matrix;
\item $\mathcal M^N_{\lambda,\Lambda}:=\left\{A\in\mathds S^N:\eig_A\subset[\lambda,\Lambda]^N\right\}$.
\end{itemize}

\section{Proof of the Result}

We start defining what we mean with \emph{viscosity solution} to the second order PDE problem
\begin{equation}\label{eq:defpde}
F\left(x,u(x),\nabla u(x),D^2u(x)\right)=0,
\end{equation}
where $F$ is a arbitrary continuous function. For a detailed overview of the viscosity solution theory we refer to \cite{userguide}.

\begin{defin}
Given an upper semicontinuous function $u$ we say that a function $\varphi$ is a \emph{supertangent} to $u$ at $x$ if $x$ is a local maximizer of $u-\varphi$.\\
Similarly we say that a function $\psi$ is a \emph{subtangent} to a lower semicontinuous function $v$ at $x$ if $x$ is a local minimizer of $v-\psi$.
\end{defin}

\begin{defin}
An upper semicontinuous function $u$ is called a \emph{viscosity subsolution} to \eqref{eq:defpde} if, for any suitable $x$ and $C^2$ supertangent $\varphi$ to $u$ at $x$,
\[
F\left(x,u(x),\nabla\varphi(x),D^2\varphi(x)\right)\ge0.
\]
Similarly a lower semicontinuous function $v$ is called a \emph{viscosity supersolution} to \eqref{eq:defpde} if, for any suitable $x$ and $C^2$ subtangent $\psi$ to $v$ at $x$,
\[
F\left(x,v(x),\nabla\psi(x),D^2\psi(x)\right)\le0.
\]
Finally a continuous function $u$ is called a \emph{viscosity solution} to \eqref{eq:defpde} if it is both a super and a subsolution to \eqref{eq:defpde}.
\end{defin}

The next definition is crucial to our problem, as we will show later.

\begin{defin}
We will call a set $C\subset\mathds R^N$ a \emph{convex cone} if for every $x,y\in C$ then $x+y\in C$ and $\alpha x\in C$ for any non negative $\alpha$.\\
We will say that a set $D$ satisfies the \emph{exterior cone condition} if, for any $x\in\partial D$, there is a convex cone $C$ with $\interior C\neq\emptyset$ and a positive $\delta$ such that $(x+C)\cap\overline D\cap B_\delta(x)=\{x\}$.
\end{defin}

We will deal with the following problem:

\begin{prob}\label{pprob}
Let $D\subset\mathds R^N$ be an open bounded set which satisfies an exterior cone condition. Given $\lambda$, $\Lambda$ in $\mathbb R$ such that $0<\lambda\le\Lambda$, we consider the Pucci's extremal operator $\mathcal P^+_{\lambda,\Lambda}$ defined by
\[
\mathcal P^+_{\lambda,\Lambda}(S):=\max_{A\in\mathcal M_{\lambda,\Lambda}^N}\langle A,S\rangle,\qquad\text{for any }S\in\mathds S^N.
\]
We further consider two continuous functions $f,g:\mathds R^N\to\mathds R$, and denote by $\ell$ a positive constant such that $\max\limits_{x\in\overline D}(|f(x)|\vee|g(x)|)\le\ell$.\\
We want to study the solution $u$ to the elliptic PDE
\[
\left\{
\begin{aligned}
&\dfrac12\mathcal P^+_{\lambda,\Lambda}\left(D^2u\right)+f(x)=0,&\qquad&x\in D,\\
&u(x)=g(x),&&x\in\partial D.
\end{aligned}
\right.
\]
\end{prob}

For the above problem a comparison result holds true, as a consequence of the maximum principle proved in \cite{CaCa95}.

\begin{theo}\label{pcomp}
Let $u$ and $v$ be respectively a subsolution and a supersolution to \cref{pprob} such that $u\le v$ on $\partial D$. Then $u\le v$ on $\overline D$.
\end{theo}

Usually, to obtain representation formulas for viscosity solutions to a second order PDE with linear operator, is useful to use a matrix $\sigma$ such that $\sigma\sigma^\dag$ is the diffusion part of the operator. We point out that every positive semidefinite matrix can be decomposed in this way. Using a similar approach, we define the set $K$ made up by the matrices $\sigma\in\mathds R^{N\times N}$ such that $\sigma\sigma^\dag\in\mathcal M_{\lambda,\Lambda}^N$, so that we can write
\begin{equation}\label{eq:puccimaxaux}
\mathcal P^+_{\lambda,\Lambda}(S):=\max_{\sigma\in K}\left\langle\sigma\sigma^\dag,S\right\rangle.
\end{equation}
For each $\sigma\in K$, $t\in[0,\infty)$ and $\zeta\in L^2\left(\Omega,\mathcal F_t;\mathds R^N\right)$, we define the linear operator $L_\sigma$ such that, for any $S\in\mathds S^N$,
\[
L_\sigma(S):=\frac12\left\langle\sigma\sigma^\dag,S\right\rangle,
\]
the stochastic process
\begin{equation}\label{eq:sde}
X_{\sigma,s}^{t,\zeta}:=\zeta+\int_t^s\sigma dW_r,\qquad s\in[t,\infty),
\end{equation}
and $\tau_\sigma^{t,\zeta}$ as the exit time of $X_\sigma^{t,\zeta}$ from $D$, i.e.,
\begin{equation}\label{eq:exittime}
\tau_\sigma^{t,\zeta}:=\inf\left\{s\in[0,\infty):X_{\sigma,t+s}^{t,\zeta}\notin D\right\}.
\end{equation}
For notation's sake we will omit the dependence from $\sigma$ of $X$ and $\tau$ when obvious and we will write $X_{\sigma,s}^\zeta$ and $\tau_{\sigma,s}^\zeta$ instead of $X_{\sigma,s}^{0,\zeta}$ and $\tau_{\sigma,s}^{0,\zeta}$, respectively.

It is well known from the Feynman--Kac formula that the viscosity solution to the problem
\[
\left\{
\begin{aligned}
&L_\sigma\left(D^2v\right)+f(x)=0,&\qquad&x\in D,\\
&v(x)=g(x),&&x\in\partial D,
\end{aligned}
\right.
\]
is given by
\[
v(x):=\mathds E\left(g\left(X^x_{\sigma,\tau_\sigma^x}\right)+\int_0^{\tau_\sigma^x}f\left(X_{\sigma,t}^x\right)dt\right).
\]
Moreover, since $\dfrac12\mathcal P^+_{\lambda,\Lambda}=\max\limits_{\sigma\in K}L_\sigma$, basic properties of viscosity solutions yield that
\[
u(x):=\sup_{\sigma\in K}\mathds E\left(g\left(X^x_{\sigma,\tau_\sigma^x}\right)+\int_0^{\tau_\sigma^x}f\left(X_{\sigma,t}^x\right)dt\right)
\]
is a subsolution to \cref{pprob}.

Our method to obtain representation formulas relies on a \emph{dynamic programming principle} which is an adaptation of the one presented by Denis, Hu and Peng in \cite[Proposition 45]{art:pathpeng} and is based on a construction on a broader set which contains $K$. This set, which we call $\mathcal A$, is made up of the progressive processes $\sigma:[0,\infty)\times\Omega\to\mathds R^{N\times N}$ which are cadlag, i.e. right continuous and left bounded, on $[0,\infty)$ and such that, for any $t\in[0,\infty)$ and $\omega\in\Omega$, the eigenvalues of $(\sigma\sigma^\dag)(t,\omega)$ belong to $[\lambda,\Lambda]^N$. $\mathcal A$ is obviously non empty, since it contains $K$. Furthermore we have by our assumptions that, for any $p>0$ and $T\in[0,\infty)$, $\mathcal A\subset L^p([0,T]\times\Omega)$, thus we endow $\mathcal A$ with with the topology of the $L^2$--convergence on compact set, which is to say that a sequence in $\mathcal A$ converges to an element of $\mathcal A$ if and only if it converges in $L^2([0,T]\times\Omega)$ for any $T\in[0,\infty)$.\\
For any a.e. finite stopping time $\rho$, an useful subset of $\mathcal A$, which we will use later, is $\mathcal A^\rho$, which consists of the $\sigma$ belonging to $\mathcal A$ such that $\{\sigma_{\rho+t}\}_{t\in[0,\infty)}$ is progressive with respect to the filtration $\left\{\mathcal F_t^\rho\right\}_{t\in[0,\infty)}$. Trivially $\mathcal A^0=\mathcal A$. We point out that if $\sigma\in\mathcal A^\rho$, then the process $\left\{X_{\sigma,\rho+t}^{\rho,x}\right\}_{t\in[0,\infty)}$ is progressive with respect to the filtration $\left\{\mathcal F_t^\rho\right\}_{t\in[0,\infty)}$. Moreover, as a consequence of the definition, we have
\[
\mathcal P^+_{\lambda,\Lambda}(S)=\max_{\sigma\in\mathcal A^\rho}\left\langle\sigma_{\rho+t}\sigma^\dag_{\rho+t},S\right\rangle,
\]
since this is true for each $\omega\in\Omega$ and $t\in[0,\infty)$, thanks to \eqref{eq:puccimaxaux}.

The dynamic programming principle is, in our case as in \cite{art:pathpeng}, an instrument that permit us to break a stochastic trajectory in two or more parts (this intuitively explain why we require $\sigma$ to be cadlag in time), i.e., for any a.e. finite stopping time $\rho$,
\begin{equation}\label{eq:pdpp}
\begin{aligned}
u(x)=&\sup_{\sigma\in\mathcal A}\mathds E\left(Y_\sigma^\rho\left(X^x_{\sigma,\rho\wedge\tau_\sigma^x}\right)+\int_0^{\rho\wedge\tau_\sigma^x}f\left(X_{\sigma,s}^x\right)ds\right)\\
=&\sup_{\sigma\in\mathcal A}\sup_{\sigma'\in\mathcal A}\mathds E\left(Y^\rho_{\sigma'}\left(X^x_{\sigma,\rho\wedge\tau_\sigma^x}\right)+\int_0^{\rho\wedge\tau_\sigma^x}f\left(X_{\sigma,s}^x\right)ds\right)\\
=&\sup_{\sigma\in\mathcal A}\mathds E\left(\esss_{\sigma'\in\mathcal A}\mathds E\left(Y^\rho_{\sigma'}(y)\middle|y=X^x_{\sigma,\rho\wedge\tau_\sigma^x}\right)+\int_0^{\rho\wedge\tau_\sigma^x}\!f\left(X_{\sigma,s}^x\right)ds\right)\!,
\end{aligned}
\end{equation}
where
\[
Y^\rho_\sigma(\zeta):=g\left(X^{\rho,\zeta}_{\sigma,\rho+\tau_\sigma^{\rho,\zeta}}\right)+\int_0^{\tau_\sigma^{\rho,\zeta}}f\left(X_{\sigma,\rho+s}^{\rho,\zeta}\right)ds.
\]
To prove it, we proceed by steps analyzing, for any $\zeta\in L^2\left(\Omega,\mathcal F_\rho;\overline D\right)$, the functions $Y^\rho_\sigma(\zeta)$ and $\Phi_\rho(\zeta):=\esss\limits_{\sigma\in\mathcal A}\mathds E(Y^\rho_\sigma(\zeta)|\mathcal F_\rho)$. We point out that the initial datum $\zeta$ here represents the first part of the trajectory defined by $X_\sigma$ broken off at $\rho$, i.e. is a generalization of the term $X_{\sigma,\rho}^x$ in \eqref{eq:pdpp}, moreover
\begin{equation}\label{eq:ptimeun}
\tau_\sigma^{\rho,X_{\sigma,\rho}^x}=\max\{\tau_\sigma^x-\rho,0\}.
\end{equation}
Further notice that, given two a.e. finite stopping times $\rho$ and $\rho'$ such that $\rho\le\rho'$, the definition yields
\begin{equation}\label{eq:concY}
Y^\rho_\sigma(\zeta)=Y^{\rho'}_\sigma\left(X^{\rho,\zeta}_{\sigma,\rho'\wedge\tau_\sigma^{\rho,\zeta}}\right)+\int_0^{\rho'\wedge\tau_\sigma^{\rho,\zeta}}f\left(X_{\sigma,\rho+s}^{\rho,\zeta}\right)ds.
\end{equation}

\begin{rem}
The reason why we require the elements of $\mathcal A$ to be cadlag in time can be understood looking at the second identity of \eqref{eq:pdpp}. While it is clearly true in $\mathcal A$, since given two elements $\sigma_1$, $\sigma_2$ in $\mathcal A$ $\sigma_{3,t}:=\sigma_{1,t}\chi_{\{t<\rho\}}+\sigma_{2,t}\chi_{\{t\ge\rho\}}$ belongs to $\mathcal A$, this could not be valid in a space with time continuous elements.
\end{rem}

We preliminarily study the continuity of $Y_\sigma^\rho$ and $\tau_\sigma^\rho$, where $\rho$ is an a.e. finite stopping time.

\begin{lem}\label{ptimecontaux}
For any a.e. finite stopping time $\rho$, $x\in\mathds R^N$ and $\sigma\in\mathcal A$, define the stopping time $
\overline\tau^{\rho,x}_\sigma:=\inf\left\{s\in[0,\infty):X^{\rho,x}_{\sigma,\rho+s}\notin\overline D\right\}$. Then
\begin{equation}\label{eq:ptimecontaux.1}
\mathds P(\tau^{\rho,x}=\overline\tau^{\rho,x})=1
\end{equation}
for any a.e. finite stopping time $\rho$, $x\in\mathds R^N$ and $\sigma\in\mathcal A$.
\end{lem}

The identity \eqref{eq:ptimecontaux.1} will be used in \cref{sdetimecont,sdetimecontX} to prove the continuity of the exit times and consequently will allow us to prove the continuity of our candidate viscosity solution. In the proof of the Lemma we will employ the exterior cone condition for $D$.

\proof
The statement is obvious if $x\notin\overline D$, while, if instead $x\in\overline D$, by \eqref{eq:ptimeun} it is equivalent to
\begin{equation}\label{eq:ptimecontaux1}
\mathds P\left(\overline\tau^{\tau^{\rho,x},X^{\rho,x}_{\tau^{\rho,x}}}=0\right)=1,
\end{equation}
for any a.e. finite stopping time $\rho$ and $\sigma\in\mathcal A$. Instead of \eqref{eq:ptimecontaux1} we will prove the stronger result
\begin{equation}\label{eq:ptimecontaux2}
\mathds P(\overline\tau^{\rho,y}=0)=1
\end{equation}
for any a.e. finite stopping time $\rho$, $\sigma\in\mathcal A$, $x\in\overline D$ and $y\in\partial D$.\\
We will proceed by steps.\\
\emph{Step 1.} Using the same proof of \cite[Proposition III.3.1]{basspde} we have that, for any $\sigma\in\mathcal A$ and $x\in\partial D$, $\mathds P(\overline\tau^x=0)=1$.\\
\emph{Step 2.} Since for any a.e. finite stopping time $\rho$ and $\sigma\in\mathcal A^\rho$ we can take a $\overline \sigma\in\mathcal A$ such that, for any $y\in\partial D$, $X^{\rho,y}_{\sigma,\rho+t}$ and $X_{\overline\sigma,t}^y$ have the same distribution, we have that $\overline\tau^{\rho,y}_\sigma$ and $\overline\tau^y_{\overline\sigma}$ have the same distribution, and consequently \eqref{eq:ptimecontaux2} is true for any a.e. finite stopping time $\rho$, $\sigma\in\mathcal A^\rho$ and $y\in\partial D$.\\
\emph{Step 3.} Fix an a.e. finite stopping time $\rho$ and consider the set
\[
\mathcal J:=\left\{
\begin{aligned}
\sigma\in\mathcal A:\sigma|_{[\rho,\infty)}=\sum_{i=0}^n\chi_{A_i}\sigma_i|_{[\rho,\infty)},\text{ where }\{\sigma_i\}_{i=0}^n\subset\mathcal A^\rho\\
\text{and }\{A_i\}_{i=0}^n\text{ is a }\mathcal F_\rho\text{--partition of }\Omega
\end{aligned}
\right\}.
\]
For each $\sigma:=\sum\limits_{i=0}^n\chi_{A_i}\sigma_i\in\mathcal J$, \eqref{eq:ptimecontaux2} holds true for any $y\in\partial D$ since, by the previous step,
\begin{align*}
\mathds P\left(\overline\tau^{\rho,y}_\sigma=0\right)=&\sum_{i=0}^n\mathds P\left(A_i\cap\left\{\overline\tau^{\rho,y}_{\sigma_i}=0\right\}\right)=\sum_{i=0}^n\mathds P(A_i)\mathds P\left(\overline\tau^{\rho,y}_{\sigma_i}=0\right)\\
=&\sum_{i=0}^n\mathds P(A_i)=1.
\end{align*}
\emph{Step 4.} As a consequence of the density of the simple functions and since each set in $\mathcal F_\infty$ is the result of intersections and unions of sets in $\mathcal F_\rho$ and $\mathcal F_\infty^\rho$, fixed a $\sigma\in\mathcal A$ we know that it is the limit of a sequence in $\mathcal J$. Therefore, fixed a $\overline\sigma\in\mathcal A$ and two positive constants $\alpha$ and $\varepsilon$, we will prove that there exists a $\sigma\in\mathcal J$ such that
\begin{equation}\label{eq:ptimecontaux3}
\mathds P(\overline\tau^{\rho,y}_{\overline\sigma}>\overline\tau^{\rho,y}_\sigma+\alpha)<\varepsilon.
\end{equation}
Notice that by the previous step $\overline\tau^{\rho,y}_\sigma=0$ a.e. for any $\sigma\in\mathcal J$, thus \eqref{eq:ptimecontaux3} is equivalent to
\[
\mathds P(\overline\tau^{\rho,y}_{\overline\sigma}>\alpha)<\varepsilon,
\]
therefore the arbitrariness of $\overline\sigma$, $y$, $\rho$, $\alpha$ and $\varepsilon$ proves \eqref{eq:ptimecontaux2} for any a.e. finite stopping time $\rho$, $\sigma\in\mathcal A$ and $y\in\partial D$ concluding the proof.\\
Fix an a.e. finite stopping time $\rho$, $y\in\partial D$ and define for any $\sigma\in\mathcal J$ the stopping times
\[
\tau^\beta_\sigma:=\inf\left\{t\in[0,\infty):\inf_{z\in D}\left|X^{\rho,y}_{\rho+t}-z\right|\ge\beta\right\}.
\]
By \cref{sdetimebound} we can take a positive $T$ depending only on $D$, $\ell$, $\lambda$ and $\varepsilon$ such that $\mathds P\left(\overline\tau^{\rho,y}_{\overline\sigma}\ge T\right)<\dfrac\varepsilon3$. Similarly we can choose a $\beta$, depending on $\alpha$ and $\varepsilon$, such that $\mathds P\left(\tau^\beta_\sigma>\overline\tau^{\rho,y}_\sigma+\alpha\right)<\dfrac\varepsilon3$ for any $\sigma\in\mathcal J$, in fact if that would not be true we should have, thanks to the reverse Fatou's lemma,
\[
\mathds P(\overline\tau^{\rho,y}_\sigma>\alpha)\ge\limsup_{\beta\to0}\mathds P\left(\tau^\beta_\sigma>\overline\tau^{\rho,y}_\sigma+\alpha\right)\ge\frac\varepsilon3
\]
for some $y\in\overline D$, in contradiction with the previous step. Thus we have that
\begin{align*}
\mathds P(\overline\tau^{\rho,y}_{\overline\sigma}>\overline\tau^{\rho,y}_\sigma+\alpha)\le&\mathds P\left(\tau^\beta_\sigma>\overline\tau^{\rho,y}_\sigma+\alpha\right)+\mathds P(\overline\tau^{\rho,y}_{\overline\sigma}\ge T)\\
&+\mathds P\left(\left\{\overline\tau^{\rho,y}_{\overline\sigma}>\tau^\beta_\sigma\right\}\cup\{\overline\tau^{\rho,y}_{\overline\sigma}<T\}\right)\\
\le&\mathds P\left(\left\{\overline\tau^{\rho,y}_{\overline\sigma}>\tau^\beta_\sigma\right\}\cup\{\overline\tau^{\rho,y}_{\overline\sigma}<T\}\right)+\frac{2\varepsilon}3.
\end{align*}
Finally Markov's and Burkholder--Davis--Gundy's inequalities yield
\begin{multline*}
\mathds P\left(\left\{\overline\tau^{\rho,y}_{\overline\sigma}>\tau^\beta_\sigma\right\}\cup\{\overline\tau^{\rho,y}_{\overline\sigma}<T\}\right)\\
\begin{aligned}
\le&\mathds P\left(\left\{\left|X^{\rho,y}_{\overline\sigma,\rho+\tau^\beta_\sigma}-X^{\rho,y}_{\sigma,\rho+\tau^\beta_\sigma}\right|\ge\beta\right\}\cap\{\overline\tau^{\rho,y}_{\overline\sigma}<T\}\right)\\
\le&\frac1{\beta^2}\mathds E\left(\sup\limits_{t\in[0,T]}\left|X^{\rho,y}_{\overline\sigma,\rho+t}-X^{\rho,y}_{\sigma,\rho+t}\right|^2\right)\le\frac c{\beta^2}\mathds E\left(\int_0^T|\overline\sigma_t-\sigma_t|^2dt\right),
\end{aligned}
\end{multline*}
where $c$ is an independent constant, hence we can choose a $\sigma\in\mathcal J$ such that $\mathds P\left(\left\{\overline\tau^{\rho,y}_{\overline\sigma}>\tau^\beta_\sigma\right\}\cup\{\overline\tau^{\rho,y}_{\overline\sigma}<T\}\right)<\dfrac\varepsilon3$ proving \eqref{eq:ptimecontaux3}.
\endproof

The previous lemma permit us to use \cref{sdetimecont,sdetimecontX} to prove that:

\begin{prop}\label{ptimecont}
The function $(t,x,\sigma)\in[0,\infty)\times\mathds R^N\times\mathcal A\mapsto\tau_\sigma^{t,x}$ is, under our assumptions, continuous in probability.
\end{prop}

\begin{lem}\label{pviscontaux}
Under our assumptions there exists a constant $c$, which depends only on $\ell$, $\lambda$ and $D$ such that
\begin{equation}\label{eq:pvisbound}
\mathds E\left(|Y^\rho_\sigma(\zeta)|^2\right)\le c
\end{equation}
for any a.e. finite stopping time $\rho$, $\sigma\in\mathcal A$ and $\zeta\in L^2\left(\Omega,\mathcal F_\rho;\overline D\right)$. Furthermore the function
\begin{align*}
Y:[0,\infty)\times\mathds R^N\times\mathcal A&\longrightarrow L^1(\Omega;\mathds R)\\
(t,x,\sigma)&\longmapsto Y^t_\sigma(x)
\end{align*}
is continuous and, for any a.e. finite stopping time $\rho$, the family of functions $\{Y^\rho_\sigma\}_{\sigma\in\mathcal A}$ is uniformly equicontinuous.
\end{lem}
\proof
\Cref{sdetimebound} yields \eqref{eq:pvisbound}, while the continuity is just a consequence of \eqref{eq:pvisbound}, \cref{ptimecont} and the dominated convergence theorem. Thus we only have to show the equicontinuity.\\
To prove the uniform equicontinuity we fix a positive $\varepsilon$ and find a $\delta$ for which
\[
\mathds E(|Y^\rho_\sigma(x)-Y^\rho_\sigma(y)|)<\varepsilon
\]
for any $\sigma\in\mathcal A$, $x\in\mathds R^N$ and $y\in B_\delta(x)$. Preliminarily fix an a.e. finite stopping time $\rho$ and note that, by Hölder's and Markov's inequalities and \cref{sdetimebound}, we can choose a $T>0$ such that
\[
\mathds E(|Y^\rho_\sigma(x)-Y^\rho_\sigma(y)|;\tau^{\rho,x}\ge T)=\left(\mathds E\left(|Y^\rho_\sigma(x)-Y^\rho_\sigma(y)|^2\right)\right)^\frac12(\mathds P(\tau^{\rho,x}\ge T))^\frac12\!<\frac\varepsilon3
\]
for any $x,y\in\mathds R^N$ and $\sigma\in\mathcal A$. Similarly we can, given an $\alpha>0$ that will be fixed later, use \cref{sdetimecont} to obtain a $\delta$ such that, for any $\sigma\in\mathcal A$, $x\in\mathds R^N$ and $y\in B_\delta(x)$,
\[
\mathds E(|Y^\rho_\sigma(x)-Y^\rho_\sigma(y)|;|\tau^{\rho,x}-\tau^{\rho,y}|>\alpha)<\frac\varepsilon3.
\]
Therefore we just have to prove that
\[
\mathds E(|Y^\rho_\sigma(x)-Y^\rho_\sigma(y)|;\{|\tau^{\rho,x}-\tau^{\rho,y}|\le\alpha\}\cap\{\tau^{\rho,x}<T\})<\frac\varepsilon3.
\]
Now define $\tau_x:=\tau^{\rho,x}\wedge T$, $\tau_y:=(\tau_x-\alpha)\vee(\tau^{\rho,y}\wedge(\tau_x+\alpha))$ and let
\begin{gather*}
Y_x:=g\left(X^{\rho,x}_{\rho+\tau_x}\right)+\int_0^{\tau_x}f\left(X^{\rho,x}_{\rho+s}\right)ds,\\
Y_y:=g\left(X^{\rho,y}_{\rho+\tau_y}\right)+\int_0^{\tau_y}f\left(X^{\rho,y}_{\rho+s}\right)ds.
\end{gather*}
Since $Y_x=Y^\rho_\sigma(x)$ and $Y_y=Y^\rho_\sigma(y)$ on $\{|\tau^{\rho,x}-\tau^{\rho,y}|\le\alpha\}\cap\{\tau^{\rho,x}<T\}$, we have that
\[
\mathds E(|Y^\rho_\sigma(x)-Y^\rho_\sigma(y)|;\{|\tau^{\rho,x}-\tau^{\rho,y}|\le\alpha\}\cap\{\tau^{\rho,x}<T\})\le\mathds E(|Y_x-Y_y|),
\]
so to conclude the proof we will show that, for a suitable $\delta$,
\begin{equation}\label{eq:pviscontaux2}
\mathds E(|Y_x-Y_y|)<\frac\varepsilon3
\end{equation}
for any $\sigma\in\mathcal A$, $x\in\mathds R^N$ and $y\in B_\delta(x)$.\\
For any fixed $\sigma\in\mathcal A$, from our assumptions it follows that
\begin{equation}\label{eq:pviscontaux3}
\begin{aligned}
\mathds E(|Y_x-Y_y|)\le\mathds E&\left(\vphantom{\int_0^{\tau_x\wedge\tau_y}}\left|g\left(X^{\rho,x}_{\rho+\tau_x}\right)-g\left(X^{\rho,y}_{\rho+\tau_y}\right)\right|+\alpha\ell\right.\\
&\quad+\left.\int_0^{\tau_x\wedge\tau_y}\left|f\left(X^{\rho,x}_{\rho+t}\right)-f\left(X^{\rho,y}_{\rho+t}\right)\right|dt\right),
\end{aligned}
\end{equation}
hence to prove \eqref{eq:pviscontaux2} we will give an upper bound to the right side of this inequality. Notice that
\begin{align*}
\mathds E\left(\left|X^{\rho,x}_{\rho+\tau_x}-X^{\rho,y}_{\rho+\tau_y}\right|^2\right)\le&2\mathds E\left(\left|X^{\rho,x}_{\rho+\tau_x}-X^{\rho,y}_{\rho+\tau_x}\right|^2+\left|X^{\rho,y}_{\rho+\tau_x}-X^{\rho,y}_{\rho+\tau_y}\right|^2\right)\\
\le&2\mathds E\left(|x-y|^2+\int_{(\tau_x-\alpha)^+}^{\tau_x+\alpha}|\sigma_t|^2dt\right)\\
\le&2|\delta|^2+4N\Lambda\alpha,
\end{align*}
thus there exists a $\delta'$ depending on $\alpha$ and $\delta$ such that
\begin{equation}\label{eq:pviscontaux4}
\mathds E\left(\left|g\left(X^{\rho,x}_{\rho+\tau_x}\right)-g\left(X^{\rho,y}_{\rho+\tau_y}\right)\right|;\left|X^{\rho,x}_{\rho+\tau_x}-X^{\rho,y}_{\rho+\tau_y}\right|>\delta'\right)<\frac\varepsilon6.
\end{equation}
From the Heine--Cantor theorem we know that $f$ and $g$ are uniformly continuous in $\overline D$, then, thanks to \eqref{eq:pviscontaux4}, we can choose $\delta$ and $\alpha$ such that \eqref{eq:pviscontaux3} holds true.
\endproof

We can now focus on our dynamic programming principle, which we point out is a generalization of the one presented by Denis, Hu and Peng in \cite[Section 3.1]{art:pathpeng}, and will be proved along the same lines. However our controls are cadlag in time, instead of just progressive as in \cite{art:pathpeng}, because they are more suitable to our needs, as will appear clear later in the proof of \cref{pvissol}.

Define, for any a.e. finite stopping time $\rho$ and $\zeta\in L^2\left(\Omega,\mathcal F_\rho;\mathds R^N\right)$, our candidate viscosity solution
\[
\Phi_\rho(\zeta):=\esss_{\sigma\in\mathcal A}\mathds E(Y^\rho_\sigma(\zeta)|\mathcal F_\rho).
\]
We proceed stating some preliminary results which can be proved slightly adapting the proofs of \cite[Lemmas 41--44]{art:pathpeng}.

\begin{lem}
For each $\sigma_1$ and $\sigma_2$ in $\mathcal A$ there exists a $\sigma\in\mathcal A$ such that
\[
\mathds E(Y^\rho_\sigma(\zeta)|\mathcal F_\rho)=\mathds E\left(Y^\rho_{\sigma_1}(\zeta)\middle|\mathcal F_\rho\right)\vee\mathds E\left(Y^\rho_{\sigma_2}(\zeta)\middle|\mathcal F_\rho\right).
\]
Therefore exists a sequence $\{\sigma_i\}_{i\in\mathds N}$ in $\mathcal A$ such that a.e.
\[
\mathds E\left(Y^\rho_{\sigma_i}(\zeta)\middle|\mathcal F_\rho\right)\uparrow\Phi_\rho(\zeta).
\]
We also have
\begin{equation}\label{eq:plattice.3}
\mathds E(\Phi_\rho(\zeta))=\sup_{\sigma\in\mathcal A}\mathds E(Y^\rho_\sigma(\zeta))<\infty,
\end{equation}
and, for any stopping time $\rho'$ such that $\rho'\le\rho$,
\begin{equation}\label{eq:plattice.4}
\mathds E\left(\esss_{\sigma\in\mathcal A}\mathds E(Y^\rho_\sigma(\zeta)|\mathcal F_\rho)\middle|\mathcal F_{\rho'}\right)=\esss_{\sigma\in\mathcal A}\mathds E(Y^\rho_\sigma(\zeta)|\mathcal F_{\rho'}).
\end{equation}
\end{lem}

\begin{lem}\label{pPhidet}
For each $x\in\mathds R^N$, $u(x):=\Phi_\rho(x)$ is a deterministic function. Furthermore $\Phi_\rho(x)=\Phi_0(x)$, i.e. the definition of $u$ does not depend on $\rho$.
\end{lem}

At a first sight the fact that $\Phi_\rho(x)$ is deterministic may seem trivial, but we remember to the reader that in general the conditional expectation, hence $\mathds E\left(Y^\rho_\sigma(x)\middle|\mathcal F_\rho\right)$, is not deterministic.

Thanks to \eqref{eq:plattice.3} and the equicontinuity proved in \cref{pviscontaux} the following \namecref{pviscont} holds true:

\begin{prop}\label{pviscont}
The function $u(x):=\Phi_\rho(x)=\sup\limits_{\sigma\in\mathcal A}\mathds E(Y^\rho_\sigma(x))$ is bounded and continuous.
\end{prop}

\begin{lem}\label{pPhi=u}
For each $\zeta\in L^2\left(\Omega,\mathcal F_\rho;\overline D\right)$, we have that $u(\zeta)=\Phi_\rho(\zeta)$ a.e..
\end{lem}

Now we have all the ingredients ready to prove the dynamic programming principle, however, before that, we will show another consequence of the above preliminary results.

\begin{prop}\label{pvissub}
The function $u(x):=\sup\limits_{\sigma\in\mathcal A}\mathds E\left(Y^0_\sigma(x)\right)$ is a continuous viscosity subsolution to \cref{pprob}.
\end{prop}
\proof
Thanks to \cref{pviscont} we already know that $u$ is continuous. In order to demonstrate that $u$ is a viscosity subsolution, we will prove that, for any fixed a $\sigma\in\mathcal A$, $v(t,x):=\mathds E(Y_\sigma^t(x))$ is a viscosity subsolution to
\begin{equation}\label{eq:pvissub1}
\partial_tv+\frac12\mathcal P^+_{\lambda,\Lambda}\left(D_x^2v\right)+f=0
\end{equation}
on $\{0\}\times\mathds R^N$. Indeed, thanks to well known properties of the viscosity solution theory, $u(t,x):=\sup\limits_{\sigma\in\mathcal A}\mathds E\left(Y_\sigma^t(x)\right)$ is then a viscosity subsolution to \eqref{eq:pvissub1} on $\{0\}\times\mathds R^N$, while \eqref{eq:plattice.3} and \cref{pPhidet} yield that $\partial_tu\equiv0$. This implies that $u$ is a viscosity subsolution to \cref{pprob}.\\
We argue by contradiction assuming that, given an $x\in D$ and a $\varphi$ supertangent to $v$ in $(0,x)$,
\[
\partial_t\varphi(0,x)+\frac12\left\langle\sigma_0\sigma_0^\dag,D^2_x\varphi(0,x)\right\rangle\le\partial_t\varphi(0,x)+\frac12\mathcal P^+_{\lambda,\Lambda}\left(D^2_x\varphi(0,x)\right)<-f(x).
\]
We assume that $\varphi(0,x)=v(0,x)$ and, for an opportune $\delta>0$,
\begin{equation}\label{eq:pvissub2}
\varphi(t,y)\ge v(t,y),\qquad\text{for any }(t,y)\in[0,\delta)\times B_\delta(x)\subseteq D.
\end{equation}
Since $\sigma$ is cadlag, there exists a set $A\in\mathcal F$ with $\mathds P(A)>0$ such that
\begin{equation}\label{eq:pvissub3}
\partial_t\varphi(t,y)+\frac12\left\langle\sigma_t\sigma^\dag_t,D_x^2\varphi(t,y)\right\rangle<-f(y)
\end{equation}
a.e. on $A$ for any $(t,y)\in[0,\delta)\times B_\delta(x)$, possibly taking a smaller $\delta>0$. Let $\rho$ be the stopping time
\[
\rho:=\delta\wedge\inf\{t\in[0,\infty):|X^x_t-x|\ge\delta\}
\]
and note that by \eqref{eq:concY}
\begin{equation}\label{eq:pvissub4}
\begin{aligned}
v(0,x)=\mathds E\left(\!Y^\rho_\sigma\left(X^x_\rho\right)+\!\int_0^\rho\!f(X_t^x)dt\right)\!=\mathds E\left(\!v\left(\rho,X^x_\rho\right)+\!\int_0^\rho\!f(X_t^x)dt\right)\!.
\end{aligned}
\end{equation}
From Itô's formula we obtain
\[
\varphi(0,x)=\mathds E\left(\varphi\left(\rho,X^x_\rho\right)-\int_0^\rho\left(\frac12\left\langle\sigma_t\sigma^\dag_t,D^2_x\varphi(t,X_t^x)\right\rangle+\partial_t\varphi(t,X_t^x)\right)dt\right),
\]
thus, by \eqref{eq:pvissub2}, \eqref{eq:pvissub3} and \eqref{eq:pvissub4}, $v(0,x)<\varphi(0,x)$ for any $\sigma\in\mathcal A$, in contradiction with the dynamic programming principle.
\endproof

Now we bring back our focus to the dynamic programming principle.

\begin{theo}[Dynamic Programming Principle]\label[dpp]{pdpp}
Let $\rho$ be an a.e. finite stopping time, then under our assumptions, for any $x\in\mathds R^N$,
\begin{equation}\label{eq:pdpp.1}
u(x)=\sup_{\sigma\in\mathcal A}\mathds E\left(u\left(X^x_{\sigma,\rho\wedge\tau_\sigma^x}\right)+\int_0^{\rho\wedge\tau_\sigma^x}f\left(X_{\sigma,t}^x\right)dt\right).
\end{equation}
\end{theo}
\proof
By \eqref{eq:concY} and since the $\sigma$ in $\mathcal A$ are cadlag, we have
\[
u(x)=\sup_{\sigma\in\mathcal A}\mathds E\left(Y_\sigma^0(x)\right)=\sup_{\sigma\in\mathcal A}\sup_{\sigma'\in\mathcal A}\mathds E\left(Y^\rho_{\sigma'}\left(X^x_{\sigma,\rho\wedge\tau_\sigma^x}\right)+\int_0^{\rho\wedge\tau_\sigma^x}f\left(X_{\sigma,t}^x\right)dt\right).
\]
Furthermore it follows from \eqref{eq:plattice.4} and \cref{pPhi=u} that
\begin{align*}
\sup_{\sigma'\in\mathcal A}\mathds E\left(Y^\rho_{\sigma'}\left(X^x_{\sigma,\rho\wedge\tau_\sigma^x}\right)\right)=&\mathds E\left(\esss_{\sigma'\in\mathcal A}\mathds E\left(Y^\rho_{\sigma'}\left(X^x_{\sigma,\rho\wedge\tau_\sigma^x}\right)\middle|\mathcal F_\rho\right)\right)\\
=&\mathds E\left(u\left(X^x_{\sigma,\rho\wedge\tau_\sigma^x}\right)\right),
\end{align*}
hence \eqref{eq:pdpp.1} is true.
\endproof

Finally we can prove the main statement of this paper.

\begin{theo}\label{pvissol}
The function $u(x):=\sup\limits_{\sigma\in\mathcal A}\mathds E\left(Y^0_\sigma(x)\right)$ is the unique viscosity solution to \cref{pprob} such that $u(x)=g(x)$ for any $x\in\partial D$.
\end{theo}
\proof
That $u(x)=g(x)$ on $\partial D$ follows from the definition, while the uniqueness is a consequence of \cref{pcomp}. From \cref{pvissub} we already know that $u$ is a continuous viscosity subsolution, hence we only have to show that $u$ satisfies the supersolution property.\\
Fixed $x\in D$, let $\psi$ be a subtangent to $u$ in $x$ which we assume, without loss of generality, equal to $u$ in $x$ and $\delta$ a positive constant such that
\begin{equation}\label{eq:psolvis1}
\psi(y)\le u(y),\qquad\text{for any }y\in B_\delta(x)\subseteq D.
\end{equation}
We will proceed by contradiction assuming that
\[
\frac12\mathcal P^+_{\lambda,\Lambda}\left(D^2\psi(x)\right)>-f(x).
\]
We know by definition that there exists a continuous and deterministic $\sigma\in\mathcal A$ for which
\[
\mathcal P^+_{\lambda,\Lambda}\left(D^2\psi(x)\right)=\left\langle\sigma_0\sigma^\dag_0,D^2\psi(x)\right\rangle,
\]
then, by continuity,
\begin{equation}\label{eq:psolvis2}
\frac12\left\langle\sigma_t\sigma^\dag_t,D^2\psi(y)\right\rangle>-f(y)
\end{equation}
for any $(t,y)\in[0,\delta)\times B_\delta(x)$, possibly taking a smaller $\delta$.\\
If we define
\[
\overline Y:=u\left(X^x_\rho\right)+\int_0^\rho f(X^x_t)dt,
\]
where $\rho$ is the stopping time
\[
\rho:=\delta\wedge\inf\{t\in[0,\infty):|X^x_t-x|\ge\delta\},
\]
then we get from the \cref{pdpp} that
\begin{equation}\label{eq:psolvis3}
\sup_{\sigma\in\mathcal A}\mathds E\left(\overline Y\right)=u(x)=\psi(x).
\end{equation}
However Itô's formula yields that
\[
\psi(x)=\mathds E\left(\psi\left(X^x_\rho\right)-\frac12\int_0^\rho\left\langle\sigma_t\sigma^\dag_t,D^2\psi(X_t^x)\right\rangle dt\right),
\]
which contradicts \eqref{eq:psolvis3}, since $u\left(X^x_\rho\right)\ge\psi\left(X^x_\rho\right)$, which follows from \eqref{eq:psolvis1}, and \eqref{eq:psolvis2} imply that $\mathds E\left(\overline Y\right)>\psi(x)$.
\endproof

\begin{appendices}

\section{Exit Times}

In this appendix we will give some properties of the exit times defined in \eqref{eq:exittime} which are needed in the paper. The following result is well known and can be found in many textbook, like \cite{basspde,pinskybook}.

\begin{prop}\label{sdetimebound}
Fixed $\sigma\in\mathcal A$, let $D$ be a bounded set, $X$ as in \eqref{eq:sde} and $\tau$, as in \eqref{eq:exittime}, the exit time of $X$ from $D$. Then there exists a constant $c$, which depends only on $D$, $\ell$ and $\lambda$, such that $\mathds E\left(\tau^{\rho,x}\right)\le c$ for any a.e. finite stopping time $\rho$ and $x\in\mathds R^N$.\\
In particular $\tau^{\rho,x}$ is a.e. finite and, for any $\varepsilon>0$, there exists a $T\in[0,\infty)$, which depends only on $D$, $\lambda$ and $\varepsilon$, such that $\mathds P\left(\tau^{\rho,x}\ge T\right)<\varepsilon$ for any a.e. finite stopping time $\rho$ and $x\in\mathds R^N$.
\end{prop}

We proceed giving the following

\begin{theo}\label{sdetimecont}
Fixed $\sigma\in\mathcal A$, let $D$ be a bounded set, $X$ as in \eqref{eq:exittime}, $\tau$ the exit time of $X$ from $D$ and $\overline\tau$ the exit time of $X$ from $\overline D$. Assume that, for some a.e. finite stopping time $\rho$, $\mathds P\left(\tau^{\rho,x}=\overline\tau^{\rho,x}\right)=1$ for any $x\in\mathds R^N$. Then, fixed two positive constants $\varepsilon$ and $\alpha$ we have that exists a $\delta>0$, depending on $D$, $\lambda$, $\varepsilon$ and $\alpha$, such that
\[
\mathds P\left(\left|\tau^{\rho,x}-\tau^{\rho,y}\right|>\alpha\right)<\varepsilon
\]
for any $x\in\mathds R^N$ and $y\in B_\delta(x)$. Therefore $\tau^\rho$ is continuous in probability with respect to $x$.
\end{theo}

The continuity of the exit times is already known, see, e.g., \cite{basspde}. What we add here is a sort of uniformity with respect to the initial conditions.

\proof
To ease notation assume that $\rho=0$ and restrict $x$ in $\overline D$, since for $x\notin\overline D$ this is obviously true. We start fixing $x\in\overline D$, defining for any $y\in\overline D$
\[
\tau_\beta^y:=\inf\left\{t\in[0,\infty):\inf_{z\in D}|X^y_t-z|\ge\beta\right\}
\]
and noting that
\begin{align*}
\mathds P(|\tau^x-\tau^y|>\alpha)=&\mathds P(\{|\tau^x-\tau^y|>\alpha\}\cap\{\tau^x\ge T\})\\
&+\mathds P\left(\{\tau^y>\tau^x+\alpha\}\cap\left\{\tau^x_\beta>\tau^x+\alpha\right\}\cap\{\tau^x<T\}\right)\\
&+\mathds P\left(\{\tau^y>\tau^x+\alpha\}\cap\left\{\tau^x_\beta\le\tau^x+\alpha\right\}\cap\{\tau^x<T\}\right)\\
&+\mathds P\left(\{\tau^x>\tau^y+\alpha\}\cap\left\{\tau^y_\beta>\tau^y+\alpha\right\}\cap\{\tau^x<T\}\right)\\
&+\mathds P\left(\{\tau^x>\tau^y+\alpha\}\cap\left\{\tau^y_\beta\le\tau^y+\alpha\right\}\cap\{\tau^x<T\}\right)\\
\le&\mathds P(\tau^x\ge T)+\mathds P\left(\tau^x_\beta>\tau^x+\alpha\right)\\
&+\mathds P\left(\left\{\tau^y>\tau^x_\beta\right\}\cap\left\{\tau^x_\beta<T+\alpha\right\}\right)+\mathds P\left(\tau^y_\beta>\tau^y+\alpha\right)\\
&+\mathds P\left(\left\{\tau^x>\tau^y_\beta\right\}\cap\left\{\tau^y_\beta<T\right\}\right).
\end{align*}
By \cref{sdetimebound} we can take a positive $T$ depending only on $D$, $\ell$, $\lambda$ and $\varepsilon$ such that $\mathds P(\tau^x\ge T)<\dfrac\varepsilon5$. Similarly we can choose a $\beta$, depending on $\alpha$ and $\varepsilon$, such that $\mathds P\left(\tau^y_\beta>\tau^y+\alpha\right)<\dfrac\varepsilon5$ for any $y\in\overline D$, in fact if that would not be true we should have, thanks to the reverse Fatou's lemma,
\[
\mathds P(\overline\tau^y>\tau^y+\alpha)\ge\limsup_{\beta\to0}\mathds P\left(\tau^y_\beta>\tau^y+\alpha\right)\ge\frac\varepsilon5
\]
for some $y\in\overline D$, in contradiction with our hypothesis.\\
Now, for the other terms, we can use Markov's inequality to get
\begin{align*}
\mathds P\left(\left\{\tau^y>\tau^x_\beta\right\}\cap\left\{\tau^x_\beta<T+\alpha\right\}\right)\le&\mathds P\left(\left\{\left|X^x_{\tau^x_\beta}-X^y_{\tau^x_\beta}\right|\ge\beta\right\}\cap\left\{\tau^x_\beta<T+\alpha\right\}\right)\\
\le&\frac1{\beta^2}\mathds E\left(\sup\limits_{t\in[0,T+\alpha]}\left|X^x_t-X^y_t\right|^2\right)\\
\le&\frac1{\beta^2}|x-y|^2
\end{align*}
and similarly
\[
\mathds P\left(\left\{\tau^x>\tau^y_\beta\right\}\cap\left\{\tau^y_\beta<T\right\}\right)\le\frac1{\beta^2}|x-y|^2.
\]
Therefore there exists a $\delta>0$ depending on $D$, $\lambda$, $\varepsilon$ and $\alpha$ such that
\[
\mathds P\left(\left\{\tau^y>\tau^x_\beta\right\}\cap\left\{\tau^x_\beta<T+\alpha\right\}\right)+\mathds P\left(\left\{\tau^x>\tau^y_\beta\right\}\cap\left\{\tau^y_\beta<T\right\}\right)\le\dfrac{2\varepsilon}5
\]
for any $y\in B_\delta(x)$ and consequently $\mathds P\left(\left|\tau^{t,x}-\tau^{t,y}\right|>\alpha\right)<\varepsilon$.
\endproof

\begin{rem}\label{sdetimecontX}
We point out that under the same assumptions and with a similar proof, we can prove that the function $(t,\sigma)\in[0,\infty)\times\mathcal A\mapsto\tau_\sigma^{t,x}$ is continuous for any $x\in\mathds R^N$.
\end{rem}

\end{appendices}

\phantomsection

\pdfbookmark[1]{References}{References}
\bibliography{biblio}
\end{document}